# CONTRACTIBLE OPEN 3-MANIFOLDS WHICH NON-TRIVIALLY COVER ONLY NON-COMPACT 3-MANIFOLDS

## Robert Myers

### 1. Introduction

Suppose $M$ is a closed, connected, orientable, irreducible 3-manifold such that $G = \pi_1(M)$ is infinite. One consequence of Thurston's geometrization conjecture is that the universal covering space $\widetilde{M}$ of $M$ must be homeomorphic to $\mathbf{R}^3$. This has been verified directly under several different additional assumptions on $G$. (See, for example, [2], [3], [6], [19].) Since $\widetilde{M}$ is irreducible [11] and contractible and $G$ is a torsion free group which acts as a group of covering translations on $\widetilde{M}$, one alternative approach to proving that $\widetilde{M}$ is homeomorphic to $\mathbf{R}^3$ has been to study such group actions on Whitehead manifolds (irreducible, contractible open 3-manifolds which are not homeomorphic to $\mathbf{R}^3$) in an attempt to show that they cannot cover compact 3-manifolds. The author showed in [14] that genus one Whitehead manifolds (those which are monotone unions of solid tori) cannot admit any non-trivial such action. Wright generalized this in [21] to the class of eventually end-irreducible Whitehead manifolds, a class which includes all those which are monotone unions of cubes with a bounded number of handles. Tinsley and Wright [18] then gave specific examples of non-eventually end-irreducible Whitehead manifolds which admit no such actions. They also gave a counting argument which shows that there must exist Whitehead manifolds which are infinite cyclic covering spaces of other non-compact 3-manifolds but cannot cover compact 3-manifolds. However, their argument does not provide any *specific* such examples. In this paper we present a general method for analyzing torsion free groups $G$ of covering translations on certain Whitehead manifolds $W$ and use it to provide specific examples of this type.

The method is based on a result of Geoghegan and Mihalik [4] which implies that there is an isomorphic image of $G$ in the mapping class group $\mathcal{H}(W)$ of $W$. The examples are constructed so that every torsion free subgroup of $\mathcal{H}(W)$ which is isomorphic to the fundamental group of a closed, irreducible 3-manifold $M$ must be a group for which it is known that $\widetilde{M}$ must be homeomorphic to $\mathbf{R}^3$.

The paper is organized as follows. Section 2 sets up the application of the Geoghegan-Mihalik theorem to the problem. Section 3 describes the general structure of the examples and gives three properties which together imply that the

---

Typeset by $\mathcal{A}_{\mathcal{M}}\mathcal{S}$-TEX

mapping class group has the appropriate structure. The examples are all obtained by taking a contractible non-compact 3-manifold $V$ whose boundary consists of two disjoint planes, gluing these planes together to get a 3-manifold with infinite cyclic fundamental group, and then letting $W$ be its universal covering space. It is thus a union of copies $V_i$ of $V$ stacked end to end such that $V_i \cap V_{i+1}$ is a plane $E_i$. The first property is that none of the $E_i$ are "trivial" in $W$ and no two distinct $E_i$ are parallel in $W$. The second property is that every element of $\mathcal{H}(W)$ has a representative which preserves this family of planes. The third property is that $\mathcal{H}(V)$ is a torsion group with a bound on the orders of its elements. Section 4 quotes results of [16] and [17] to give conditions on $V$ which imply the first two properties. Section 5 describes conditions on an exhaustion of $V$ which imply the third property. Section 6 gives specific examples of such exhaustions.

## 2. The General Method

We refer to [7] or [8] for general definitions in 3-manifold topology.

Given an orientable manifold $X$ we let $Homeo^+(X)$ denote the group of orientation preserving homeomorphisms of $X$. The mapping class group $\mathcal{H}(X)$ is then the quotient of $Homeo^+(X)$ by the normal subgroup consisting of those homeomorphisms which are isotopic to the identity. We let $\Phi : Homeo^+(X) \to \mathcal{H}(X)$ be the quotient epimorphism. If $Y$ is a submanifold of $X$ we let $\mathcal{H}(X \, rel \, Y)$ denote the group of orientation preserving homeomorphisms of $X$ which are the identity on $Y$ modulo isotopies which are the identity on $Y$.

We let $Fr_X Y$ and $Int \, _X Y$ denote the topological boundary and interior of $Y$ in $X$, with the subscripts deleted when they are clear from the context. The manifold theoretic boundary and interior of $Y$ are denoted by $\partial Y$ and $int \, Y$.

An exhaustion $\{K_n\}$ for a connected, non-compact manifold $W$ is a sequence of compact, connected, codimension zero submanifolds of $W$, indexed by the non-negative integers, whose union is $W$ such that $K_n \subseteq Int \, K_{n+1}$, $K_n \cap \partial W$ is either empty or a codimension zero submanifold of $\partial W$, and $W - Int \, K_n$ has no compact components.

A map $f : M \to N$ between manifolds is **end-proper** if pre-images of compact sets are compact. It is $\partial$-**proper** if $f^{-1}(\partial N) = \partial M$. It is **proper** if it has both these properties. These terms are applied to submanifolds if their inclusion maps have the given property.

**Theorem 2.1 (Geoghegan-Mihalik).** *Let $W$ be a Whitehead manifold. Suppose $G$ is a torsion free group of covering translations of $W$. Then the restriction of $\Phi$ to $G$ is one to one.*

*Proof.* Define a space $X$ to be **strongly connected at** $\infty$ if any two end-proper rays in $X$ are end-proper homotopic. (A ray in $X$ is a map from $[0, \infty)$ to $X$.) Theorem B of [4] states that if $G$ is a torsion free group of covering translations of a connected manifold $X$ such that for every $g \in G$ the quotient $X/ <g>$ is non-compact, and $X$ is not strongly connected at $\infty$, then the restriction of $\Phi$ to $G$ is one to one.

In our case $W/<g>$ is non-compact because otherwise it would be homeomorphic to $S^1 \times S^2$ (see Theorem 5.2 of [7]), thereby contradicting the asphericity of $W$. We claim that $W$ is not strongly connected at $\infty$. Let $\{K_n\}$ be an exhaustion for $W$. Choose an end-proper ray in $W$ for which the image of $[n, \infty)$ lies in $W - K_n$ for all $n$ and use inclusion and change of basepoint along it to define the sequence of homomorphisms $\pi_1(W-K_0) \xleftarrow{\varphi_0} \pi_1(W-K_1) \xleftarrow{\varphi_1} \pi_1(W-K_2) \xleftarrow{\varphi_2} \cdots$. If $W$ is strongly connected at $\infty$, then by Proposition 1.1 of [4] this sequence is semi-stable, i.e. for every $m$ there exists an $n \geq m$ such that for every $k \geq n$ the image of $G_k \to G_m$ is equal to the image of $G_n \to G_m$. This implies that $W$ is end 1-movable in the sense of Brin and Thickstun [1] and thus by Corollary (A) to their Theorem 1.1 $W$ is homeomorphic to $\mathbf{R}^3$. □

We now state a sufficient criterion for a universal cover to be $\mathbf{R}^3$ which we shall use in the examples.

**Theorem 2.2.** *Let $M$ be a closed, orientable, irreducible 3-manifold with $\pi_1(M)$ infinite. If $\pi_1(M)$ contains a subgroup of finite index which either has infinite abelianization or contains a non-trivial normal abelian subgroup, then the universal covering space of $M$ is homeomorphic to $\mathbf{R}^3$.*

*Proof.* By passing to a finite sheeted covering space we may assume that the subgroup is $\pi_1(M)$. If the abelianization is infinite then $M$ is Haken, and so the result follows from Waldhausen [19]. So assume that $A$ is a non-trivial normal abelian subgroup of $\pi_1(M)$. If $A$ contains a $\mathbf{Z} \oplus \mathbf{Z}$ subgroup then we are done by Hass-Rubinstein-Scott [6]. If $A$ is $\mathbf{Z}$, then we are done by the proof of the Seifert fibered space conjecture [2] or [3].

If neither of these hold then $A$ is a subgroup of $\mathbf{Q}$. (See Theorem 9.14 of [7].) Conjugation by elements of $\pi_1(M)$ gives a homomorphism $\psi : \pi_1(M) \to Aut(A)$. Let $I$ be the ideal $\mathbf{Z} \cap A$ in $\mathbf{Z}$. Choose a generator $m$ for $I$. Then for $\alpha \in Aut(A)$, $\alpha(m) = p/q$. If $k \in I$, then $k = rm$, and so $\alpha(k) = r\alpha(m) = rp/q = k(p/mq)$. Then for $k/\ell \in A$, $\ell\alpha(k/\ell) = \alpha(k)$, and so $\alpha(k\ell) = (k/\ell)(p/mq)$. This shows that $Aut(A)$ is a subgroup of $\mathbf{Q} - \{0\}$.

Let $\mathbf{Q}^+$ be the group of positive rationals, and let $K = \psi^{-1}(\mathbf{Q}^+)$. Then $K$ has index at most two in $\pi_1(M)$ and contains $A$. If $\psi(K)$ is non-trivial, then $K$ has infinite abelianization and so we are done by Waldhausen. If $\psi(K)$ is trivial, then $A$ lies in the center of $K$, and so we are done by the Seifert fibered space conjecture. □

## 3. The Construction and Main Theorem

**Construction 3.1.** *Let $V$ be an irreducible, contractible, non-compact 3-manifold such that $\partial V$ consists of two disjoint planes $E$ and $E'$. Let $V^*$ be the 3-manifold obtained by identifying $E$ with $E'$ by an orientation reversing homeomorphism. Let $E^*$ be the image of $E$ in $V^*$. Let $p : W \to V^*$ be the universal covering space of $V^*$. Let $\{E_i\}$ be the set of components of $E = p^{-1}(E^*)$. Let $\{V_i\}$ be the set of copies of $V$ into which $E$ splits $W$. Index so that $V_i \cap V_{i+1} = E_i$.*

Let $X$ be an open 3-manifold. An end-proper plane $P$ in $X$ is **trivial** if some component of $X - P$ has closure homeomorphic to $\mathbf{R}^2 \times [0, \infty)$, with $P = \mathbf{R}^2 \times \{0\}$. Two disjoint end-proper planes $P_0$ and $P_1$ in $X$ are **parallel** if some component of $X - (P_0 \cup P_1)$ has closure homeomorphic to $\mathbf{R}^2 \times [0, 1]$, with $P_i = \mathbf{R}^2 \times \{i\}$.

**Theorem 3.2.** *Let $W$ be as in Construction 3.1. Assume that:*
  (1) *No $E_i$ is trivial, and no distinct $E_i$ and $E_j$ are parallel.*
  (2) *Every $g \in Homeo^+(W)$ is isotopic to an $h$ such that $h(E) = E$.*
  (3) *$\mathcal{H}(V)$ is a torsion group with a bound on the orders of its elements.*

*Then every torsion free subgroup of $\mathcal{H}(W)$ has a subgroup of finite index which either has an infinite abelianization or a non-trivial normal abelian subgroup.*

**Corollary 3.3.** *Let $W$ be as in Theorem 3.2. Then $W$ cannot cover a compact 3-manifold.* □

*Proof of Theorem 3.2.* Let $D_\infty = <x, y : x^2 = 1, xyx = y^{-1}>$ be the infinite dihedral group. Regard it as the group of automorphisms of the simplicial complex $\Gamma$ obtained by triangulating $\mathbf{R}$ using $\mathbf{Z}$ as the set of vertices. We associate $i \in \mathbf{Z}$ with $V_i$ and $[i, i+1]$ with $E_i$. Define $\rho : \mathcal{H}(W) \to D_\infty$ as follows. Given $g \in Homeo^+(W)$, let $\rho(\Phi(g))$ be the automorphism of $\Gamma$ determined by the isotopic homeomorphism $h$ provided by (2). One sees that this is well defined as follows. By Theorem 5 of [20] disjoint, ambient isotopic, non-trivial, end-proper planes in an irreducible 3-manifold are parallel. Hence by (1) distinct $E_i$ and $E_j$ are not ambient isotopic. If $k$ were a homeomorphism isotopic to $g$ which determined a different automorphism of $\Gamma$, then $k$ and $h$ would send some $E_i$ to different planes; hence these planes would be ambient isotopic.

There is an obvious epimorphism $\psi : \prod \mathcal{H}(V_i \, rel \, \partial V_i) \to \ker \rho$. There are also homomorphisms $\theta_i : \mathcal{H}(V_i \, rel \, \partial V_i) \to \mathcal{H}(V_i)$ obtained by allowing isotopies which need not be the identity on $\partial V_i$. These induce $\theta : \prod \mathcal{H}(V_i \, rel \, \partial V_i) \to \prod \mathcal{H}(V_i)$ with $\ker \theta = \prod \ker \theta_i$. It follows from Lemma 3.4 below that $\ker \theta$ is abelian. By (3) we have that $\prod \mathcal{H}(V_i)$ is a torsion group.

Now suppose $G$ is a torsion free subgroup of $\mathcal{H}(W)$. By passing to a subgroup of index at most two we may assume that $\rho(G)$ lies in the subgroup $<y>$ of $D_\infty$. If $\rho(G)$ is non-trivial, then $G$ has infinite abelianization, so assume $G \subseteq \ker \rho$.

The subgroup $G \cap \psi(\ker \theta)$ of $G$ is normal and abelian. Suppose it is trivial. Let $\gamma$ be a non-trivial element of $G$. Then $\gamma = \psi(\delta)$ for some $\delta \in \prod \mathcal{H}(V_i \, rel \, \partial V_i)$. Then $\theta(\delta^k) = \theta(\delta)^k = 1$ for some $k > 1$, so $\delta^k \in \ker \theta$, and so $\gamma^k = \psi(\delta)^k = \psi(\delta^k) = 1$, contradicting the fact that $G$ is torsion free. □

**Lemma 3.4.** $\ker \theta_i$ *is abelian.*

*Proof.* Let $g$ be a homeomorphism representing an element of $\ker \theta$. Let $C^* = \partial V_i \times [0, 2]$ be a collar on $\partial V_i$ in $V_i$, with $\partial V_i \times \{0\} = \partial V_i$. By the isotopy uniqueness of collars we may assume that $g$ is the identity on $C^*$. Let $C = \partial V_i \times [0, 1]$. Let $V_i'$ be the closure of $V_i - C$, and let $V_i^*$ be the closure of $V_i - C^*$.

We claim that $g|_{V_i'}$ is isotopic to the identity. Define $k : V_i' \to V_i$ to be the identity on $V_i^*$ and to be $k(x, s) = (x, 2s - 2)$ on $\partial V_i \times [1, 2]$. Let $g_t$ be an isotopy with $g_0$ the identity of $V_i$ and $g_1 = g$. Then $k^{-1} g_t k$ is the required isotopy $h_t$ with $h_0$ the identity of $V_i'$ and $h_1 = g|_{V_i'}$.

We next claim that $g$ is isotopic to a homeomorphism $f$ which is the identity outside $C$. Define $f_t$ to be $h_t$ on $V_i'$, while on $C$ we have $f_t(x, s) = (ph_{(1-s)+st}(x, 1), s)$, where $p$ is projection onto the first factor. Then $f_1 = g$ and $f_0$ is the desired $f$.

Thus $\ker \theta$ is the image of $\mathcal{H}(C \, rel \, \partial C)$ in $\mathcal{H}(V_i \, rel \, \partial V_i)$. The result then follows from the result below. □

**Lemma 3.5.** $\mathcal{H}(\mathbf{R}^2 \times [0, 1] \, rel \, \mathbf{R}^2 \times \{0, 1\})$ *is infinite cyclic.*

*Proof.* Let $f$ be a homeomorphism of $\mathbf{R}^2 \times [0, 1]$ which is the identity on the boundary. Let $\{B_n\}$ be an exhaustion of $\mathbf{R}^2$ by concentric disks. Let $C_n = B_n \times [0, 1]$ and $F_n = (\partial B_n) \times [0, 1]$.

We may assume that $f(C_0) \subseteq Int \, C_1$, $C_1 \subseteq Int \, f(C_2)$, and $f(C_2) \subseteq Int \, C_3$. Since $C_3 - Int \, C_1$ is a product $S^1 \times [0, 1] \times [0, 1]$ with $S^1 \times [0, 1] \times \{0\} = F_1$ and $S^1 \times [0, 1] \times \{1\} = F_3$ and $f(F_2)$ is incompressible in this manifold it follows (say from Corollary 3.2 of [19]) that there is an ambient isotopy in $C_3 - Int \, C_1$, fixed on the boundary, which carries $f(F_2)$ to $F_2$. We may repeat this argument with $C_{3k}$, $C_{3k+1}$, $C_{3k+2}$, and $C_{3k+3}$ for $k \geq 1$ to get an isotopy of $f$ rel $\mathbf{R}^2 \times \{0, 1\}$ taking $f$ to a homeomorphism which carries each $C_{3k+2}$ to itself.

Thus we may assume that $f(C_n) = C_n$ for all $n$. The restriction of $f$ to $F_n$ is then isotopic rel $\partial F_n$ to a Dehn twist $\sigma(\varphi, s) = (\varphi + 2\pi k_n s, s)$, where $\varphi$ is the angular variable in $S^1$, $s \in [0, 1]$, and $k_n$ is an integer. Consider the meridional disk $\{0\} \times [0, 1] \times [0, 1]$ of the solid torus $C_{n+1} - Int \, C_n = S^1 \times [0, 1] \times [0, 1]$. It follows from the fact that $F$ must carry this disk to another meridional disk that we must have all $k_n$ equal to a constant $k$. The restriction of $f$ to the boundary determines the isotopy class relative to the boundary of the restriction of $f$ to $C_0$ and each $C_{n+1} - Int \, C_n$. It follows that $f$ is isotopic rel $\mathbf{R}^2 \times \{0, 1\}$ to the homeomorphism $\sigma(\varphi, s) = (\varphi + 2\pi k, s)$. Thus $\mathcal{H}(\mathbf{R}^2 \times [0, 1] \, rel \, \mathbf{R}^2 \times \{0, 1\})$ is cyclic.

Suppose $\sigma$ is isotopic to the identity rel $\mathbf{R}^2 \times \{0, 1\}$. Identify $\mathbf{R}^2 \times \{0\}$ with $\mathbf{R}^2 \times \{1\}$ to obtain $\mathbf{R}^2 \times S^1$. Then the map $\tau$ induced by $\sigma$ is isotopic to the identity. Let $K_n$ and $T_n$ be the images of $C_n$ and $F_n$ in $\mathbf{R}^2 \times S^1$. We may assume that the track of $T_1$ under the isotopy misses $K_0$. Let $\{\alpha, \beta\}$ be a basis for $H_1(T_1)$, where $\alpha$ is represented by an $S^1$ fiber and $\beta$ bounds in $K_1$. Then $i_*(\alpha) = i_*(\alpha + k\beta)$, where $i : T_1 \to (\mathbf{R}^2 \times S^1) - K_0$ is the inclusion. Since $i_*$ is an isomorphism we have $k = 0$. Thus $\mathcal{H}(\mathbf{R}^2 \times [0, 1] \, rel \, \mathbf{R}^2 \times \{0, 1\})$ is infinite cyclic. □

## 4. Planes in $W$

In this section we assemble results from [16] and [17] to show how to obtain manifolds $V$ for which Construction 3.1 produces a 3-manifold $W$ with properties (1) and (2) of Theorem 3.2.

Let $V$ be a non-compact 3-manifold such that $\partial V$ is either empty or a disjoint union of planes. $V$ is **aplanar** if every proper plane $P$ in $V$ is either trivial or $\partial$-**parallel** (cobounds an end-proper $\mathbf{R}^2 \times [0, 1]$ with a component of $\partial V$). A **partial** **plane** is a non-compact, simply connected surface with non-empty boundary. $V$ is **strongly aplanar** if it is aplanar and for every proper surface $\mathcal{P}$ in $V$ such that each component of $\mathcal{P}$ is a partial plane there is a collar on $\partial V$ which contains $\mathcal{P}$. $V$ is **anannular at** $\infty$ if for every compact subset $K$ of $V$ there is a larger compact subset $L$ of $V$ such that $V - L$ is anannular, i.e. every proper, incompressible annulus in $V - L$ is $\partial$-parallel.

A **plane sum** of oriented 3-manifolds $V_i$ whose boundaries consist of disjoint planes is a 3-manifold obtained by identifying the boundary planes in pairs via orientation reversing homeomorphisms, where the pattern of identification is determined by a countable tree $\Gamma$ whose vertices correspond to the $V_i$ and whose edges correspond to the pairs of planes. Construction 3.1 is a simple example of such a plane sum. The plane sum is **non-degenerate** if none of the summing planes are trivial or $\partial$-parallel and if distinct summing planes are not parallel. A summand $V_i$ is **bad** if there is a component $E_j$ of $\partial V_i$ such that $E_j \cup int\, V_i$ is homeomorphic to $\mathbf{R}^2 \times [0, \infty)$.

**Lemma 4.1.** *Plane sums which have no bad summands are non-degenerate.*

*Proof.* This is Corollary 3.3 of [17]. □

A plane sum is **strong** if it is non-degenerate and each summand is strongly aplanar and anannular at $\infty$.

**Lemma 4.2.** *Let $E$ be the union of the set of summing planes in a strong plane sum $W$ along a locally finite tree. Then every homeomorphism $g$ of $W$ is isotopic to a homeomorphism $h$ such that $h(E) = E$.*

*Proof.* This is Corollary 4.4 of [17]. □

It follows from Theorem 6.1 of [16] that given any connected, irreducible, one-ended open 3-manifold $U$ and integer $\nu \geq 1$, there is a non-compact 3-manifold $V$ such that $\partial V$ consists of $\nu$ disjoint planes, $V$ is strongly aplanar and anannular at $\infty$, and $int\, V$ is homeomorphic to $U$. We will apply a special case of that construction to a certain genus one Whitehead manifold $U$ which ensures that $\mathcal{H}(V)$ is a torsion group with a bound on the orders of its elements. We state the technical result from [16] that we need to establish that $V$ is strongly aplanar and anannular at $\infty$.

Suppose $V$ is a connected, orientable, irreducible, one-ended, non-compact 3-manifold whose boundary consists of a finite number $\nu \geq 1$ of disjoint planes. An exhaustion $\{C_n\}$ for $V$ is **nice** if for all $n \geq 0$ one has that $C_n \cap \partial V$ consists of a single disk in each component of $\partial V$, $X_{n+1} = C_{n+1} - Int\, C_n$ is irreducible, $\partial$-irreducible, and anannular, each component of $Fr\, C_n$ has positive genus and negative Euler characteristic, and $V - Int\, C_n$ has one component.

**Proposition 4.3.** *If $V$ has a nice exhaustion, then $V$ is strongly aplanar and anannular at $\infty$.*

*Proof.* This follows from Theorem 5.3 and Lemma 1.3(6) of [16]. □

5. The Mapping Class Group of $V$

**Construction 5.1.** Let $U$ be a genus one Whitehead manifold with an exhaustion $\{K_n\}$ by solid tori. Let $T_n = \partial K_n$. Let $Y_{n+1} = K_{n+1} - \mathrm{Int}\, K_n$. Let $\ell_n$ be an oriented non-contractible simple closed curve on $T_n$ which bounds an oriented surface in $Y_{n+1}$. Let $\alpha_{n+1}^0$ and $\alpha_{n+1}^1$ be disjoint proper arcs in $Y_{n+1}$ missing $\ell_n \cup \ell_{n+1}$ such that $\alpha_{n+1}^j$ joins $T_n$ to $T_{n+1}$, and $\alpha_{n+1}^j \cap T_{n+1} = \alpha_{n+2}^j \cap T_{n+1}$. Let $N_{n+1}^j$ be a regular neighborhood of $\alpha_{n+1}^j$ in $Y_{n+1}$, chosen so that $N_{n+1}^0 \cap N_{n+1}^1 = \emptyset$, and $N_{n+1}^j \cap T_{n+1} = N_{n+2}^j \cap T_{n+1}$. Let $X_{n+1} = Y_{n+1} - \mathrm{Int}\,(N_{n+1}^0 \cup N_{n+1}^1)$, $C_0 = K_0$, and $C_{n+1} = C_0 \cup X_1 \cup \cdots \cup X_{n+1}$. Let $V = \cup C_n$, $F_n = Fr_V C_n$, $A_{n+1}^j = Fr_U N_{n+1}^j$, $D^j = N_1^j \cap C_0$, and $E^j = D^j \cup (\cup A_{n+1}^j)$.

A compact, connected, orientable 3-manifold is **excellent** if it is irreducible, $\partial$-irreducible, anannular, and atoroidal (every proper, incompressible torus is $\partial$-parallel).

**Theorem 5.2.** Let $U$ and $V$ be as in Construction 5.1. Assume that:
  (1) $X_{n+1}$ and $Y_{n+1}$ are excellent.
  (2) There is no orientation preserving homeomorphism from $Y_1$ to $Y_2$ which takes $\ell_0$ to a curve homologous on $T_1$ to $\pm \ell_1$.
  (3) There is a partition $\mathcal{N}_1 \cup \mathcal{N}_2$ of the set of positive integers into two infinite subsets such that $1 \in \mathcal{N}_1$, $2 \in \mathcal{N}_2$, and for each $n+1 \in \mathcal{N}_i$, there is an orientation preserving homeomorphism $h_{n+1} : X_{n+1} \to X_i$ which takes $\ell_{n+1}$, $F_n$, $F_{n+1}$, $A_{n+1}^0$, and $A_{n+1}^1$ to $\ell_i$, $F_{i-1}$, $F_i$, $A_i^0$, and $A_i^1$, respectively.
  (4) If $n+1$ and $k+1$ are both in the same subset and $n+2$ and $k+2$ are both in the same subset, then $h_{k+1}^{-1} h_{n+1}$ and $h_{k+2}^{-1} h_{n+2}$ agree on $F_{n+1}$.
  (5) There do not exist integers $s \neq 0$ and $k_0 \geq 1 - s$ such that $k$ and $k+s$ are in the same subset for all $k \geq k_0$.

Then $V$ is strongly aplanar and anannular at $\infty$, and $\mathcal{H}(V)$ is a torsion group with a bound on the order of its elements.

*Proof.* Condition (1) and Proposition 4.3 imply the statement about $V$. We next prove a special case of the statement about $\mathcal{H}(V)$.

**Lemma 5.3.** There is an $L > 0$ such that if $g$ is an orientation preserving homeomorphism of $V$ for which $g(C_n) = C_n$ for all $n$, then $g^L$ is isotopic to the identity.

*Proof.* By Proposition 27.1 of [10] the mapping class group of an excellent 3-manifold is finite. Let $L = |\mathcal{H}(X_1)| \cdot |\mathcal{H}(X_2)|$. So $g^L|_{X_{n+1}}$ is isotopic to the identity. Since $\partial X_{n+1}$ has genus greater than one the isotopy can be chosen so that it preserves $F_n$, $F_{n+1}$, $A_{n+1}^0$, and $A_{n+1}^1$. Since each $F_{n+1}$ has negative Euler characteristic it follows from [5] that the isotopies can be chosen to match up on each $F_{n+1}$. Since $C_0$ is a solid torus with $\partial C_0 = F_0 \cup D^0 \cup D^1$ one can extend the isotopy over $C_0$. □

The remainder of the proof shows how to reduce to the special case. Let $M_{n,n+q} = X_n \cup X_{n+1} \cup \cdots \cup X_{n+q}$. Conditions (3) and (4) imply that for a fixed $q$ there are only finitely many homeomorphism types of such manifolds. By Lemma 2.1 of [15] each $M_{n,n+q}$ is excellent.

Let $F$ be a torus with two holes. Any collection of disjoint, non-contractible, non-$\partial$-parallel, mutually non-parallel simple closed curves on $F$ has at most two elements.

There is a result of Haken which states that given an excellent 3-manifold $M$ and an integer $k$, there are, up to isotopy, only finitely many proper, incompressible, $\partial$-incompressible surfaces $S$ in $M$ with $\chi(S) \geq k$. A proof of this theorem for the special case in which $M$ is closed is given as Corollary 2.3 of [9]. The proof of the general case is a straightforward generalization of that proof. We apply this result to conclude that there is an integer $H(M)$ such that any collection of proper, incompressible, $\partial$-incompressible, mutually non-isotopic surfaces in $M$ each of which is homeomorphic to $F$ must have at most $H(M)$ elements. Let $H = \max H(M_{n-3,n+4})$ for $n \geq 4$.

Now let $g$ be an orientation preserving homeomorphism of $V$. There exist $q > p > 0$ such that for $0 \leq t \leq H$ one has $g^t(C_0) \subseteq C_p$ and $g^{-t}(C_p) \subseteq C_q$. Consider $n \geq \max\{q, p+4\}$. Then $g^t(F_n)$ is an incompressible surface $F$ in $V - Int\, C_p$. Isotop $F$ in $V - Int\, C_p$ so that $\partial F$ misses $\partial F_k$ for all $k > p$ and among all such surfaces $F$ meets $\cup_{k>p} F_k$ in a minimal number of components.

*Claim 1.* $F \cap (X_n \cup X_{n+1}) \neq \emptyset$. If not, then one can extend $g^t$ to a homeomorphism $f$ of $U$ such that $f(K_n) \supseteq K_{n+1}$ or $f(K_n) \subseteq K_{n-1}$. Since the $Y_{k+1}$ are excellent it follows from Lemma 3.3 of [14] that $f$ can be isotoped so that $f(K_k) = K_{k+s}$ for some $s \neq 0$ and all sufficiently large $k$. In particular by (5) there is a $k+1 \in \mathcal{N}_1$ and a $w+1 \in \mathcal{N}_2$ such that $f$ takes $Y_{k+1}$ to $Y_{w+1}$ and $T_{k+1}$ to $T_{w+1}$. Then $f(\ell_{k+1})$ is homologous on $T_{w+1}$ to $\pm \ell_{w+1}$. Since $h_{k+1}$ and $h_{w+1}$ extend to homeomorphisms from $Y_{k+1}$ and $Y_{w+1}$ to $Y_1$ and $Y_2$ carrying $\ell_{k+1}$ and $\ell_{w+1}$ to $\ell_1$ and $\ell_2$, respectively, we have contradicted (2).

*Claim 2.* $F \subseteq M_{n-2,n+3}$ If not, then it must meet at least three of the surfaces $F_s$ with $n-3 \leq s \leq n+3$. Hence either some component of $F \cap (\cup_{k>p} F_k)$ is contractible or $\partial$-parallel on $F$ or two such components are parallel on $F$. The irreducibility and $\partial$-irreducibility of the $Y_{k+1}$ imply that all contractible components could have been removed by an isotopy, contradicting minimality. Ananannularity implies the same for $\partial$-parallel components. Suppose $A$ is an annulus in $F$ such that $\partial A$ consists of components of the intersection. We may assume that $A$ meets no other such components. If $\partial A \subseteq F_k$ for some $k$, then $A$ is parallel to an annulus in $F_k$, so we may assume there are no such annuli. Thus $A$ joins some $F_k$ to $F_{k+1}$. Since the components of $\partial A$ are non-separating there must be another such annulus $A'$. We may choose $A$ and $A'$ so that there is a disk with two holes $R$ in $F$ with one component each of $\partial R$ in $\partial V$, $\partial A$, and $\partial A'$ such that $int\, R$ contains no intersection curves. Since $A$ and $A'$ are $\partial$-parallel in some $X_{k+1}$ we can isotop $R$ in the adjacent $X_{w+1}$ containing it so that $\partial R$ lies in $A_{w+1}^0 \cup A_{w+1}^1$. Suppose $\partial R$ lies in $A_{w+1}^0$. It then consists of three concentric simple closed curves. Pushing the disks on $E^0$ bounded by the two inner curves slightly into $int\, V$ we obtain a proper embedding of a disk in $V$ which cobounds a 3-ball in $V$ with a disk on $E^0$. It follows that one of the curves has the opposite orientation from the other two when compared on $E^0$. Hence there is a proper torus with one hole embedded in $X_{w+1}$ with boundary a non-contractible curve in $A^0_{w+1}$. Thus there is a torus in $Y_{w+1}$ which meets $\alpha^0_{w+1}$ transversely in a single point and misses $\alpha^1_{w+1}$, which is impossible. Similar analyses apply to the other possibilities for $\partial R$. Thus there is no such annulus.

*Claim 3. F is $\partial$-incompressible in $M_{n-3,n+4}$.* Suppose $D$ is a $\partial$-compressing disk. Put $D$ in general position with respect to $F_{n-3} \cup F_{n+3}$. The irreducibility and $\partial$-irreducibility of $X_{n-3}$ and $X_{n+4}$ allow one to remove all intersections with this surface. The portion of $\partial D$ not in $F$ then lies in an annulus whose centerline is a component of $\partial F$ and so $\partial D$ can be isotoped into $F$, to which we then apply incompressibility.

We now have $H + 1$ incompressible, $\partial$-incompressible, homeomorphic surfaces in $M_{n-3,n+4}$ consisting of the isotopes of $F_n$, $g(F_n)$, ..., $g^H(F_n)$. At least two of these are isotopic. Hence $g^t(F_n)$ is isotopic to $F_n$ for some $t$, $0 \leq t \leq H$. Thus $g^{H!}(F_n)$ is isotopic to $F_n$. Note that this isotopy can be chosen to have compact support contained in $V - Int\, C_p$.

We next choose a sequence $n_0 < n_1 < n_2 < \cdots$ such that the isotopies of $g^{H!}(F_{n_k})$ to $F_{n_k}$ have disjoint supports. Thus we can perform a single isotopy which simultaneously performs this sequence of isotopies. Let $S_n$ denote the image of $g^{H!}(F_n)$ under this isotopy; thus $S_{N_k} = F_{n_k}$.

We now consider $n_k < n < n_{k+1}$. We may assume that $\partial S_n = \partial F_n$. The embedding maps of $S_n$ and $F_n$ are homotopic in $V - Int\, C_p$ by a homotopy under which the image of the boundary is in $\partial V - (\partial V \cap C_p)$. We may arrange for the homotopy to be fixed on $\partial S_n$. We claim that it may be chosen so that its image lies in $M_{n_k,n_{k+1}}$. Choose a collection of disjoint proper arcs in $S_n$ which cut it into disks. Then use the incompressibility of $F_{n_k} \cup F_{n_{k+1}}$ to modify the homotopy so that the images of these arcs miss this surface. Finally use irreducibility to extend the homotopy over the disks in $S_n$. It then follows from Corollary 5.5 of [19] that $S_n$ and $F_n$ are isotopic in $M_{n_k,n_{k+1}}$.

Continuing in this fashion one gets that $g^{H!}$ is isotopic to a homeomorphism which carries $C_n$ to itself for all sufficiently large $n$. By Lemma 5.3 we have that $g^{L \cdot H!}$ is isotopic to the identity of $V$. □

## 6. SPECIFIC EXAMPLES

**Theorem 6.1.** *There are specific examples of Whitehead manifolds each of which non-trivially covers another non-compact 3-manifold but does not cover a compact 3-manifold.*

*Proof.* Figure 1 shows a 4-component tangle $\lambda$ in a 3-ball $B$ called the true lover's 4-tangle. In [13] it is proven that the exterior of $\lambda$ is excellent. The exterior of the graph $\xi$ in Figure 2 is homeomorphic to the exterior of $\lambda$. By deleting two arcs from $\xi$ we obtain the 2-tangle $\mu$ in Figure 3 which is equivalent to the union of the second and third components of $\lambda$. It follows immediately from the proof of Proposition 4.1 of [12] that the exterior of $\mu$ is excellent.

We now identify the disks which constitute the left and right sides of the rectangular solid $B$ in order to obtain a solid torus $K$. This is done so that $\mu$ becomes a simple closed curve $\sigma$ and $\xi$ becomes a graph $\theta$ consisting of $\sigma$ together with two disjoint arcs $\alpha^0$ and $\alpha^1$ joining $\sigma$ to $\partial K$. It follows from Lemma 2.1 of [15] that the exteriors of $\sigma$ and $\theta$ in $K$ are excellent. We embed $K$ in $S^3$ in a standard way so that a line segment running along the bottom front edge of $B$ becomes a simple closed curve $\ell$ in $\partial K$ which bounds a disk in $S^3 - int\, K$.

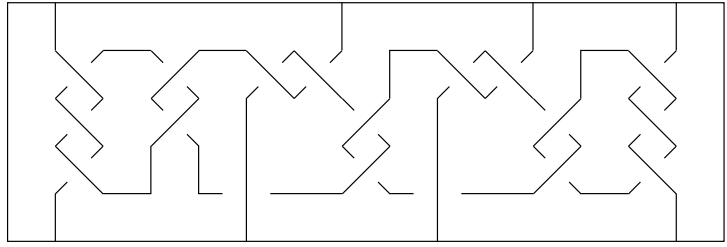

FIGURE 1. The 4-tangle $\lambda$

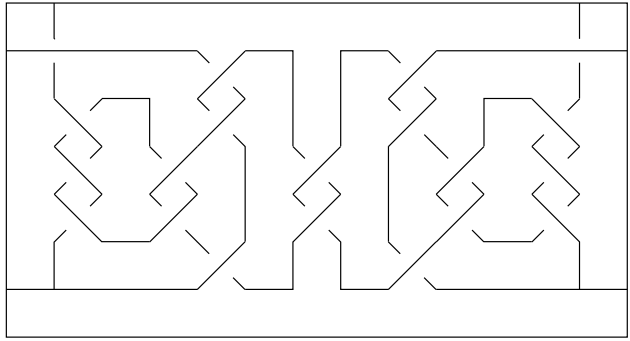

FIGURE 2. The graph $\xi$

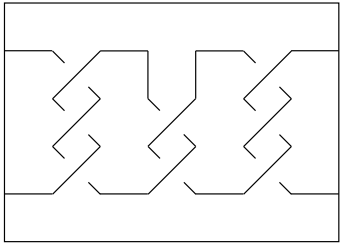

FIGURE 3. The 2-tangle $\mu$

If one changes the sense of the central clasp in the figures by changing the two overcrossings to undercrossings, thereby obtaining a new $\sigma$ and $\theta$, then the same arguments show that their exteriors in $K$ are also excellent. We denote the old and new versions by the subscripts 1 and 2, respectively.

We claim that there is no orientation preserving homeomorphism from the exterior of $\sigma_1$ to the exterior of $\sigma_2$ which carries $\ell$ to $\pm\ell$. If not, then one could extend the homeomorphism to the exteriors of these knots in $S^3$. But the two knots are the knots $8_5$ and $8_{19}$ with normalized Alexander polynomials $5 - 4(t+t^{-1}) + 3(t^2 +$ $t^{-2}) - (t^3 + t^{-3})$ and $1 - (t^2 + t^{-2}) + (t^3 + t^{-3})$, respectively, so this is impossible.

Now let $L_i$ be a regular neighborhood of $\sigma_i$ in $K$. Let $\mathcal{N}_1$ and $\mathcal{N}_2$ be as in Theorem 5.2(5), for example taking $\mathcal{N}_2$ to be the set of positive integers of the form $2^r > 1$ and $\mathcal{N}_1$ those not of this form. We construct a genus one Whitehead manifold $U$ with exhaustion $\{K_n\}$ by using as models for $(K_{n+1}, K_n)$ the pairs $(K, L_i)$ where $n+1 \in \mathcal{N}_i$. The graphs $\theta_1$ and $\theta_2$ provide the arcs $\alpha_{n+1}^j$. It is then easily checked that the construction can be carried out so as to satisfy the hypotheses of Theorem 5.2, and so the conclusion follows from Theorem 3.2. One can also vary the choice of $\mathcal{N}_1$ and $\mathcal{N}_2$ and apply Lemma 3.3 of [14] to obtain uncountably many non-homeomorphic such examples. □

Department of Mathematics, Oklahoma State University, Stillwater, OK 74078-1058, USA

*E-mail address*: myersr@math.okstate.edu